\documentclass[11pt,twoside,a4paper,notitlepage]{article}

\usepackage{graphicx}
\usepackage{mathtools}
\usepackage{enumerate}
\usepackage{amssymb}
\usepackage{amsthm}
\usepackage{mathptmx} 
\usepackage{csquotes}
\usepackage{latexsym}
\usepackage{tikz}
\usepackage{authblk}

\def\.#1{\hfill#1\kern.5em\vrule\kern-.5em}

\theoremstyle{definition}

\title{Theological Reasoning of Cantor's Set Theory}
\author{Kate\v{r}ina Trlifajov\'{a}}
\affil{       Czech Technical University in Prague\\
              Th\'{a}kurova 9, 160 00 Prague 6, Czech Republic \\
              Email:{katerina.trlifajova@fit.cvut.cz}        
}
\date{}                     

\begin{document}
\maketitle

\begin{abstract}
Discussions surrounding the nature of the infinite in mathematics have been under way for two millennia. Mathematicians, philosophers, and theologians have all taken part. The basic question has been whether the infinite exists only in potential or whether it exists in actuality. Only at the end of the 19th century a set theory was created that works with the actual infinite. Initially, this theory was rejected by other mathematicians. The creator behind the theory, the German mathematician Georg Cantor, felt all the more the need to challenge the long tradition that only recognised potential infinite. In this he received strong support from the interest among German neothomist philosophers, who, under the influence of the Encyclical of Pope Leo XIII, Aeterni Patris, began to take an interest in Cantor’s work. Gradually, his theory even acquired approval from the Vatican theologians. Cantor was able to firmly defend his work and at the turn of the 20th century he succeeded in gaining its acceptance. The storm that had accompanied its original rejection now accompanied its acceptance. The theory became the basis on which modem mathematics were and are still founded, even though the majority of mathematicians know nothing of its originally theological justification. Set theory, which today rests on an axiomatic foundation, no longer poses the question of the existence of actual infinite sets. The answer is expressed in its basic axiom: natural numbers form an infinite set. No substantiation has been discovered other than Cantor’s: the set of all natural numbers exists from eternity as an idea in God’s intellect.

\end{abstract}

\section{Prologue}
Czech students of mathematics once gave set theory the nickname \enquote{darkness},\footnote{\enquote{Darkness} is in Czech \enquote{TEMNO}, which can be an abbreviation for TEorie MNO\v zin which means a theory of sets.} and although they were only partly aware of it, they were not far from the truth. After all, the main goal of set theory is to capture and elaborate the notion of infinity in mathematics, and infinity tends to be associated with a certain darkness and ambiguity. Somewhat paradoxically, set theory seeks to remove precisely this obscurity from infinity and capture it in a clear, transparent, and static form. 

      The introduction of infinity into mathematics took more than two thousand years and involved not only mathematicians but also philosophers and theologians. It was only at the turn of the 19th and 20th centuries that the set concept, which works with the so-called actual infinity, was adopted thanks to the German mathematician Georg Cantor, marking a radical turn in the development of mathematics. After some misunderstandings, rejections, and struggles, it was accepted by the mathematical community in the early 20th century. All mathematics was built on a common set basis, which is used until today.

\section{Potential and actual infinity}

Philosophers have been concerned with the question of infinity since Antiquity, actually since geometry went beyond direct evidence and introduced logical proofs. Early on it became apparent that a rational grasp of infinity often led to paradoxes and disputes.  Aristotle formulated the fundamental problem of infinity, which eventually became central to mathematics as well as to philosophy and theology, and which has run through their history like a red thread. \emph{Does infinity exist only in possibility or also in reality?} After careful analysis, Aristotle rejected the existence of infinity in reality: \enquote{It is plain from these arguments that there is no body which is actually infinite.} (Aristotle Book III, Part 5, 206a) What remains, then, is that the infinite exists only in possibility. Similarly, in Aristotle's conception, the continuum is indeed infinitely divisible, but only potentially. 

Aristotle's question has been refined and generalized over time. Its form, with which many medieval and modern thinkers came to terms, looked like this: does infinity exist only as potential or also as actual?  A \emph{potential} infinity is one that may be ever-increasing or ever-decreasing but is never complete, whole, or finished. In contrast, the \emph{actual} or complete infinity is conceived as a whole, is given \emph{here and now} and can be treated as such. If we are dealing with infinite numbers or quantities, it is necessarily an actual infinity, since both the number and the quantity are uniquely given, complete, and determined.

    Scholastic philosophy took as its own Aristotle's radical rejection of actual infinity and expressed it in the famous thesis \emph{infinitum actu non datur}. Most medieval and modern scholars adhered to it, and until the end of the 19th century, only potential infinity was almost exclusively considered acceptable. The question of actual infinity became at the time a theological question since infinity was one of the attributes of God.  Its use in mathematics was viewed with suspicion by theologians because it might be an assault on the unique and absolute nature of God.

\section{Paradox of Reflexivity}

If one eventually accepted the actual infinity then another problem would appear. It is referred to as the paradox of reflexivity, since it concerns a one-to-one correspondence (reflection) of an actually infinite collection on its own proper part. It was demonstrated in its purest form by Galileo. (Galilei p. 40). If we compare two infinite series, the series of all natural numbers and the series of all squares of these numbers
$$1, 2, 3, \dots$$
$$1, 4, 9, \dots$$
we see on the one hand that there are fewer numbers in the second series because they are only some of the numbers in the first series. According to Euclid's Axiom 9, \enquote{the whole is greater than the part} (Euclid p. 26) the second series contains fewer numbers than the first.

On the other hand, there is the same amount of numbers in the both series. For to every number from the first series, we can always find just one its square from the second series, and vice versa, to every square from the second series we can always find just one number, its root, from the first series. According to Euclid's Axiom 8, which says: \enquote{Magnitudes which coincide with one another are equal.} (Euclid p. 26), there is the same amount of elements in both series, from this point of view they are equal. Galileo argues that the relations of being greater, less, equal, are meaningless when we speak of the infinities: \enquote{The attributes 'equal', 'greater', and 'less', are not applicable to infinite, but only to finite, quantities.} (Galilei, 41).

Leibniz pointed out a similar problem in the case of the continuum. Any two line segments $a$, $b$ contain the same number of points and yet one can be much shorter than the other. For if we form an arbitrary triangle with sides $a$, $b$ and add side $c$, then each point of side $a$ is connected by a line parallel to side $c$ with exactly one point of side $b$.

Almost until the end of the 19th century, the potential infinity was considered to be the only correct infinity, which, however, cannot be captured in mathematics otherwise than an ever-increasing or ever-decreasing finite quantity. The insufficiency of this way of expressing was apparent especially after the creation of the infinitesimal calculus, which was discovered at the end of the 17th century by Newton and Leibniz.  It treats infinitesimal quantities without introducing them exactly. The uncertainty of rules for their behaviour had sometimes led to inaccuracies in calculations and errors in proofs.

\section{Bernard Bolzano}

The only honourable exception was the Bohemian mathematician and philosopher Bernard Bolzano. Already in his works on measurable numbers from the 1930s, he treats infinity as somewhat actual. In his last book, \emph{The Paradoxes of the Infinite} (Bolzano 1851/2004), he consistently introduced and developed the mathematical theory of actual infinity. He based it on the concept of multitude [\emph{Menge}] that is a very similar sense as Cantor later used a concept of set and which has been used until today. A \emph{set} is a collection of clearly determined distinct objects into one whole.\footnote{To be precise, Bolzano's \emph{Menge} is not exactly the same notion as Cantor's \emph{Menge}. The difference is neither in its actuality nor in its being one whole but in the way the objects it contains are connected. Bolzano's \emph{multitude [Menge]} consists of parts, Cantor's \emph{set [Menge]} of elements.} In order to express infinite quantities as actual and therefore graspable, we must collect them into a whole, and this is what the meaning of a set is. A set, if it contains an infinite number of elements, is necessarily actually infinite. The key question of the existence of actual infinity is now expressed as follows: \emph{Does an infinite set exist?}

Bolzano showed that it exists in the Book 1 in \emph{Theory of Science}. (Bolzano 1937/2014, \S 30 - \S 32, p. 107 - 109). The proof is based on the notion of \emph{truth in itself} which is  Bolzano's thoroughly defined and defended concept, the cornerstone of his philosophy. He first proved \enquote{There is at least one truth in itself} by contradiction. The negation is the proposition \enquote{There are no truth} and if it was true then it contradicts itself. The proof is finished by the principle \emph{tertium non datur}.

Similarly, Bolzano proves by contradiction that there is more than one truth $A$. The proposition \enquote{There is another truth aside of $A$} must be true. And so on. Consequently there are infinitely many truths. 
Because they are truths in themselves they are valid without anybody thinking them. Moreover, God is a guarantee of this infinite multitude of truths. 
\begin{quote} Thus we must attribute to him [God] a power of knowledge that is true omniscience, that therefore comprehends an infinite multitude if truths because all truth in general. (Bolzano 1851/2006, \S 11, p. 604) \end{quote} 

Bolzano pointed out the similarity of this multitude with all natural numbers and concluded that they also form an infinite multitude as well as rational numbers, irrational numbers,  infinitely small and infinitely large quantities.  

Bolzano's theory remained almost unnoticed and without direct influence on further developments.  
Potential infinity continued to be regarded as the only acceptable infinity in mathematics.

\section{Continuum and real numbers}

There is another fundamental phenomenon that is firmly tied to infinity, although it is not obvious at first glance.  It is the phenomenon of the continuum, of continuous extension. Mathematicians and philosophers have been aware of this connection from the beginning, and knew, or at least anticipated, that once they formed one, it would necessarily be reflected in the other. The continuum and the infinite are two opposite sides of the same matter mutually connected, but the infinite is a more striking phenomenon. The basic question of infinity is this: \emph{Is the continuum infinitely divisible into continuous parts or is it composed of an infinite number of indivisible particles?}

The traditional examples of a continuum are a space, time, or motion. All of these examples can be expressed in a numerical system as a ratio to a chosen unit of a space, time, or motion. Therefore, the basic type of continuum is considered to be a straight line, which is identified with the real numbers. 
The size of any line segment can be expressed as the ratio of its length to the length of the unit line. If this ratio is rational, e.g. twice or a third, it is expressed as a given number, e.g. $2$ or $\frac{1}{3}$. However, some line segments are not in a rational ratio to the unit line, e.g. the diagonal of a square with side $1$. Their lengths correspond to irrational numbers, which are denoted symbolically, e.g. $\sqrt{2}$. The problem is how to express all these lengths arithmetically.

\section{Cantor's real numbers}

At the beginning of the 1870s, a young, talented German mathematician Georg Cantor investigated the problem of the uniqueness of trigonometric series. In doing so, he realised that a correct solution required precise definitions of irrational numbers, which at that time had not yet been established. 

He therefore expressed irrational numbers in terms of so-called \emph{fundamental sequences} of rational numbers. These are such sequences that \enquote{keep on becoming narrower}, geometrically approaching a single point. Arithmetically, a fundamental sequence is a sequence of rational numbers $(a_1, a_2, \dots , a_n, \dots)$ which has the property that for any positive rational numbers $\epsilon$ there exists an index $n$ such that for every $m > n$ it is valid $a_m - a_n < \epsilon$. Cantor assigned to such a sequence the symbolic sign $a$, which he called its limit. That this was a \enquote{symbolic} limit was soon forgotten, and the newly defined \emph{real} and old rational numbers were soon regarded as being equivalent.

The fact that such a narrowing sequence approaches geometrically only one point represented by a sign $a$ is, of course, to be stated as an axiom. Moreover, Cantor defined that two fundamental sequences have the same limit if, arithmetically speaking, their difference approaches to $0$. This step excluded infinitesimal quantities from the real numbers. Here Cantor introduced a certain interpretation of the continuum, an interpretation that eventually resulted in the introduction of the actual infinite.

\section{Sets and their cardinalities}

For an excellent mathematician such as Cantor, it became an exciting adventure to explore the newly emerging structure without prejudice. He began to work with collections of real numbers as with sets, even though he had not yet called them so and was probably not even aware of their actual infinity at the beginning. He called them \enquote{Inbegriff}, then for a long time \enquote{Mannigfaltigkeit}.

In 1882 he read Bolzano's book \emph{The Paradoxes of the Infinite} and it may encouraged him to openly admit and defend the actual infinity. He later came to call the sets \enquote{Menge}, a term used by Bolzano and common until today. In French, the present-day term \enquote{ensemble} has been used from the beginning.

The disturbing paradox of reflexivity was not considered by Cantor to be a lack, but rather a feature that characterizes infinite sets. If a \emph{one-to-one correspondence} can be established between two sets, i.e. if there exists a prescription such that every element of one set is unambiguously connected into a pair with an element from the other set, then Cantor defines that these two sets have the same size which he calls \emph{cardinality}. He accepted the one-to-one correspondence as the basis of his theory.

       He began to compare different subsets of real numbers from this point of view. Finite sets with the same number of elements have the same cardinality. The Galileo paradox demonstrates that the set of natural numbers and that of their squares have the same cardinality. Also the set of rational and that of algebraic numbers (that is, the roots of all polynomial equations, for example, the square roots) have this cardinality. This cardinality is called \emph{countable}. From the existence of a one-to-one correspondence between a countable set and the set of all natural numbers follows that the elements of a countable set can be numbered by natural numbers.

     Leibniz's example shows that there is a one-to-one correspondence between any two segments of a line, so all line segments have the same cardinality. Also the set of all real numbers has this cardinality. It is the cardinality of \emph{continuum}. Cantor proved the real numbers cannot be numbered by natural number, consequently the cardinality of continuum is not countable, it is greater. Only this theorem justifies the introduction of the notion of the cardinality for comparing sets. Otherwise, it might well be that all infinite sets would be countable. 

     Cantor proved another surprising theorem, that any $n$-dimensional Euclidean space (e.g., a square or a cube) also has the cardinality of a continuum and not a greater one. Consequently, for example, an arbitrarily small line segment contains as many points as an n-dimensional infinite space. He made commented it himself: \enquote{I see it, but I do not believe it}. He had expected that the dimension of space could be defined in terms of cardinality. 

      Cantor introduced other typical today commonly well-known set notions: intersection, union, subset. He also introduced a derivative $A'$ of a set $A$. It is the set of all its cumulation points (that is, those to which some fundamental sequence formed from points of the set $A$ converges). The second derivative $A''$ is the derivative of the first derivative, and so on. Cantor succeeded in constructing a special set $C$ which has an infinitely many  different derivatives $C, C', C'', \dots$  He denoted their intersection $C^\infty$. It was not an empty set and Cantor could make its derivative, and to it again another derivative, etc.: $C^{\infty +1}, C^{\infty +2}, \dots$  At this point Cantor realized that what he was investigated was not only the theory of real numbers, but also of actually infinite numbers. He began to devote systematically to their study.

\section{Infinite ordinal numbers}

        In 1883, Cantor published \emph{Foundations of a General Theory of Manifolds}, subtitled \emph{A Mathematical-Philosophical Study in the Theory of the Infinite} (Cantor 1883).  In the introduction, he makes a small apology for daring to depart from the long-standing traditional interpretation of infinity. On the other hand, he is otherwise unable to continue in his work. 

\begin{quote} I am so dependent on this extension of the number concept that without
it I should be unable to take the smallest step forward in the theory of sets
[Mengen]; this circumstance is the justification (or, if need be, the apology) for
the fact that I introduce seemingly exotic ideas into my work. For what is
at stake is the extension or continuation of the sequence of integers into the
infinite; and daring though this step may seem, I can nevertheless express, not
only the hope, but the firm conviction that with time this extension will have
to be regarded as thoroughly simple, proper, and natural. (Cantor 1883, p. 70).
\end{quote}

     Finite numbers are formed according to the so-called first principle of formation, that is, by adding $1$ to the previous number. Cantor laid down a second principle of formation. To an increasing sequence of integers which does not have a greatest element he created a new number which can be considered as greater than all these numbers or also as their limit. Thus he added to the series of natural numbers the smallest infinite number which is greater than all the finite numbers and denoted it by the last letter of the Greek alphabet $\omega$. 

     Using the first principle, one can then obtain from the number $\omega$ its successors $\omega + 1, \omega+2, \dots$. Using the second principle, he then got another infinite number $2\omega$, etc. The new sequence of the ordinal numbers thus formed looks like this:

$$1, 2, 3, \dots \omega, \omega+1, \omega+2, \dots  2\omega, 2\omega+1,  \dots \dots \omega \cdot \omega=\omega^2 , \dots \dots  \omega^3 , \dots \dots \omega^\omega \dots \dots$$

  Each of these numbers, which he called \emph{ordinal}, is characterised by the special type of ordering, i.e. by the way how the preceding ordinal numbers are arranged.  

The arithmetic operations of addition, multiplication and power of ordinal numbers are defined accordingly, and not surprisingly, have different properties than those operations for finite numbers. For example, neither the addition nor the multiplication are commutative. It is valid $1 + \omega = \omega$ but $\omega + 1 \neq \omega$ and $2 \cdot \omega = \omega$ but $\omega \cdot 2 \neq \omega$. On the other hand, the opposite operations of subtraction and division are not introduced. Cantor did attempt to do so, but it proved to be difficult and essentially meaningless.  

Countable sets can be arranged to sequences according to different ordinal numbers. Thus, for example, the set of all natural numbers can be arranged as follows. 

\begin{itemize}
\item According to $\omega: (1, 2, 3, \dots).$
\item According to $\omega+1 : (2, 3, 4, \dots 1).$ 
\item According to $\omega+2 : (3, 4, 5, \dots 1, 2)$ .
\item According to $2 \cdot \omega : (1, 3, 5, \dots 2, 4, 6, \dots)$.
\end{itemize}

We can see that $\omega+1, \omega+2, 2 \cdot \omega$ and other above mentioned infinite ordinal numbers are countable since they can be numbered by natural numbers which are elements of $\omega$.    

\section{Infinite cardinal numbers}

     Cantor introduced yet another kind of infinite numbers, called \emph{cardinal} numbers designating cardinalities of sets. The first smallest infinite cardinal numbers $\aleph_0$ denotes the cardinality of countable sets.

     Using a sequence of alephs numbered with ordinal numbers, Cantor marked the increasing sequence of infinite cardinalities. Each successive term denotes the first successive greater cardinality. The smallest one is the cardinal number $\aleph_0$, the successive cardinal number $\aleph_1$ is the first cardinality greater than $\aleph_0$. Then $\aleph_2$ is the first cardinality greater than $\aleph_1$, and so on. The results in an infinite sequence of increasing cardinalities.
$$\aleph_0, \aleph_1, \aleph_2, \dots \aleph_\omega, \aleph_{\omega+1}, \dots \dots \aleph_{2\omega}, \dots \dots \aleph_{\omega \cdot \omega} \dots \dots$$    

Cantor introduced arithmetic operations for cardinal numbers, again only addition, multiplication and power, which differ from these operations on finite numbers. For an infinite cardinal number $\aleph$ and a finite natural number $n$, it holds true
 $$\aleph + n = n + \aleph = \aleph, \aleph + \aleph = \aleph, n \cdot \aleph = \aleph \cdot n = \aleph.$$ 

It is relatively easy to prove that the cardinality of real numbers is $2^{\aleph_0}$, it is the cardinality of \emph{continuum}. The set of all subsets of natural numbers has the same cardinality. The important Cantor's theorem says the cardinality of continuum is greater than that countability $$2^{\aleph_0} >\aleph_0.$$  

Cantor's lifelong dream was to prove that $2^\aleph$ is $\aleph_1$, the first cardinality greater than $\aleph_0$. This conjecture is called the \emph{continuum hypotheses}. 
$$2^{\aleph_0} = \aleph_1.$$
     
\section{Cantor's justification of the existence of the actual infinity}

Cantor felt the need to confront the mathematical and philosophical tradition that had been rejected actual infinity for centuries. As an example of actual infinity, he gives his infinite ordinal numbers, the first and most important of which is $\omega$. Cantor tried to show that if previous philosophers Aristotle, Locke, Descartes, Spinoza, or Leibniz had known his new ordinal and cardinal numbers and the rules for calculating them, they would have agreed to their existence. According to Cantor, hitherto infinite numbers have been rejected on the basis of the common misconception that they must have the same properties as finite numbers. However, although they have different properties, they do not lead to a contradiction and can therefore be introduced.

      Against the scholastic principle \emph{infinitum actu non datur} Cantor puts another claim: \emph{All things, whether finite or infinite, are definite and, with the exception of God, can be determined by the intellect.} (Cantor 1883, p.76). Cantor considers the only limitation of the new mathematical ideas to be their intrinsic consistency and their consistency with already well-established notions. He praises mathematics for its relative lightness and freedom, in which it is not necessary to examine the external reality of its concepts and need only deal with their internal reality. He answers in advance a probable objection: what are infinite ordinal and cardinal numbers actually for.  For him \emph{the essence of mathematics lies in its freedom.} (Cantor, p. 79). 

      Later, however, Cantor's thinking seems to have changed, and he tried to find a different, deeper justification, more in accordance with the philosophical and theological tradition.

\section{Encyclical \emph{Aeterni Patris}} 

     If Cantor expected his revolutionary theory to arouse enthusiasm, he must have been disappointed. His results were not accepted by mathematicians. They either ignored them or rejected them outright. The main opponent, Leopold Kronecker, called Cantor's work \emph{humbug}. Through his rather behind-the-scenes dealings, he ensured that Cantor, who was a professor at the University of Halle, did not get a better position in Berlin or Gottingen. Mittag-Leffler, Cantor's previous publisher, did publish \emph{Foundations} in French, but without the philosophical passages. He advised him not to publish further results until mathematicians would accept his theory. Cantor became interested in the philosophical side of actual infinity and began to explore its ontological nature. Although actual infinity had already been mathematically established, its existence was not philosophically sufficiently justified. 

In the meantime, Cantor received unexpected support from the Roman Catholic Church. In 1879, the newly elected Pope Leo XIII issued one of his first encyclicals, \emph{Aeterni Patris}. He sought to revive the tradition of philosophical thinking, especially Thomistic philosophy. He recommended dealing with the humanities and natural sciences in the light of the teachings of St. Thomas. Under this influence, some theologians became interested in Cantor's theory. The prominent neo-Thomist theologian Constantine Gutberlet published his own work defending actual infinity. He argues that God himself provides for the existence of both irrational and infinite numbers. Since God's mind is immutable and infinite, it must contain absolute, infinite, complete sets. Either the existence of actual infinity is assumed or the infinite intellect and eternity of God's mind must be given up. 

\section{Cantor and St. Thomas Aquinas}

      The encyclical \emph{Aeterni Patris} refers especially to St. Thomas Aquinas, but he is known to have rejected actual infinity. Cantor points out that St. Thomas, in presenting five proofs of God's existence, does not use this \enquote{false thesis} at all, since he found it unreliable for this purpose.

The Thomistic doctrine was otherwise close to Cantor which is evident from many details. He used its terminology, its arguments and its strictly rational way of reasoning. Therefore, he also took St. Thomas's arguments seriously and sought to refute them. he quoted Aquinas' objection against the actual existence of an infinite multitude from the \emph{Theological Summa}. \begin{quote} But this is impossible. For every multitude must be contained under some species of multitude, and the species of multitude correspond to the species of number. But no species of number is infinite, since each number is such that it is a multitude measured by the unit. Hence, it is impossible for there to be an actual infinite multitude. (Thomas Aquinas p. 1, q. VII, a. 4) \end{quote}
    
Cantor was convinced that Thomas Aquinas and others rejected actual infinity on a basis of a fundamental error that infinite numbers would have to have the same properties as finite numbers. 
\begin{quote}{All so-called proofs against the possibility of actual infinite numbers, as can be distinctly demonstrated in every case and can also be concluded from general principles, are in the main point faulty thereby, and therein lies their \emph{protos pseudos}. (Cantor 1885).} \end{quote}
He claimed that there was only one answer to St. Thomas's objections. To introduce new infinite numbers, ordinal and cardinal, and so prove that they can exist. They even describe the world more richly than finite numbers alone.

      In the same article, Thomas Aquinas presents other arguments that Cantor rejects as unfounded and erroneous. Thomas shows first that the multitude of all numbers is potentially infinite, since we can keep multiplying and adding them among themselves and always get new numbers, or the plurality of different figures in geometrical space, since we can keep creating more and more again. However, he does raise an objection that can still be used against the actual infinity.  \begin{quote}As a day is brought into actuality not as a simultaneous whole, but instead successively; and, similarly, an infinite multitude is brought into actuality not as a simultaneous whole, but successively.\end{quote}

\section{Cantor and St. Augustine}

        Cantor appealed to St. Augustine, where he thought he found confirmation of his beliefs. St. Augustine first proves that there are infinitely many numbers, for we can always not only increase them by one, but also double or multiply them. Then he turns against those who claim that God does not know all numbers:

\begin{quote} 
Who is so left to himself as to say so?  \dots 
In our books also it is said to God, \emph{Thou hast ordered all things in number, and measure, and weight.} (Wis. 11, 21) The prophet also says, \emph{Who bringeth out their host by number.} (Is. 40, 26) And the Saviour says in the Gospel, \emph{The very hairs of your head are all numbered.} (Matt. 10, 30). Far be it, then, from us to doubt that all number is known to Him \emph{whose understanding,} according to the Psalmist, \emph{is infinite.} (Ps. 146, 5)  The infinity of number, though there be no numbering of infinite numbers, is yet not incomprehensible by Him whose understanding is infinite. And thus, if everything which is comprehended is defined or made finite by the comprehension of him who knows it, then all infinity is in some ineffable way made finite to God, for it is comprehensible by His knowledge. For his wisdom, simply manifold and uniformly manifold, 
comprehends all things with an understanding incomprehensible. (Augustine, book XII, chapter 18)\end{quote}    
     But Augustine says no more than that numbers are infinitely many and God clearly comprehends them all by an understanding incomprehensible - but he makes no mention that humans have such a capacity. Rather, he points out the difference between our and God's knowledge.

     Cantor repeatedly claims that the natural numbers, both together and separately, as well as all infinite numbers, ordinal and cardinal, exist in a higher degree of reality from eternity as thoughts in God's intellect. (Hallett 1984, p. 21), (Dauben 1990, p. 228). This argument is hard to accept today. But Cantor goes much further, attributing the ability to understand and capture the set of all natural numbers and other infinite sets to human understanding. Not only as if everything we conclude to be consistent actually \enquote{exists} in the mind of God, but also as if we ourselves could be able to look into it.

\section{Three respects of the actual infinite}

In the paper \emph{On the Various Standpoints
With Regard to the Actual Infinite} (Cantor 1886), Cantor, to avoid the danger of religious error, such as pantheism distinguishes three main respects of manifesting the actual infinite. 
\begin{enumerate}
\item in God, \emph{in natura naturans}, in the creating nature, where it is called \emph{Absolute},
\item \emph{in concreto seu in natura naturata}, in perceptible things in created nature where he calls it \emph{Transfinitum},
\item \emph{in abstracto}, that may be comprehended by human cognition in the form of  \emph{transfinite numbers}. 
\end{enumerate}
Earlier and contemporary philosophers can be divided according to which of these respects they accept and which they do not. He himself fully and unconditionally agrees with all three.
 
\begin{quote}This basis, which I consider the only right one, only a few stand; perhaps I am temporally the first, who represents this standpoint with complete determination and in all its consequences, however this I know for certain, that I shall not be the last one who defends it! \end{quote} 

 Speculative theology deals with the investigation of the \emph{Absolute} that is essentially unincreasable and therefore mathematically indeterminable. It is strictly separated from the second and third concepts which fall into metaphysics and mathematics that is to be conceived as an indeed \emph{Infinite} but nevertheless yet increasable. The confusion of the two forms frequently leads two pantheism which, for example, \enquote{constitutes the Achilles heel of Spinoza’s Ethics}. 

 An example of the actual infinite quantity in the created world is the infinite number of essences contained in the entire universe, and also the number of Leibnizian monads of which any physical body is composed.

\section{Correspondence with Cardinal Franzelin}

     Cantor's letters to the leading Vatican theologian, Cardinal Franzelin , show the importance they both attached to the proof of the existence of the actual infinite. (Cantor 1886) This problem is not a minor theological detail, but concerns the crucial question. Cantor suggests two ways in which this fact follows.

\begin{quote} One proof proceeds from the concept of God and concludes first of all from the highest Perfection of God’s Being the possibility of the creation of a Transfinitum ordinatum, then from His Benevolence and Magnificence the necessity of the actually ensued creation of a Transfinitum. \end{quote}
    
 \begin{quote} Another proof shows a posteriori, that the assumption of a Transfinitum in natura naturata renders possible a better, because more perfect explanation of the phenomena, especially the organisms and psychical manifestations, than the opposing hypothesis.
  \end{quote}   

    While Cardinal Franzelin praised Cantor for his distinction between the \emph{Absolute} and the \emph{Transfinitum}, the mention of the necessity of the creation of the \emph{Transfinitum} provoked disapproval. 
In the following letter Cantor sought to explain his position.

\begin{quote} In the brief indication of my previous letter, it was not my intention at the point in question, to speak of an objective, metaphysical necessity of the act of creation, to which God the \emph{absolute Free} would have been subjugated; on the contrary, I wanted to point to a certain subjective necessity for us, to infer from God’s Benevolence and Magnificence an actually \emph{ensued} (\emph{not} a parte Dei \emph{ensuing}) creation, not only of a \emph{Finitum ordinatum}, but also of a \emph{Transfinitum ordinatum}.  \end{quote}

This clarification was sufficient for Italian theologians to endorse Cantor's set theory. Cantor was proud of the acceptance his theory of transfinite numbers and would frequently remind the words of Cardinal Franzelin:
 \begin{quote} Thus the two concepts of  the Absolute Infinite and the Actual infinite in the created world or in the \emph{Transfinitum} are essentially different, so that in comparing the two one must only describe the former as \emph{properly infinite}, the later as equivocally infinite. When conceived in this way, so far as I can see at present, there is no danger to religious truths in your concept of \emph{Transfinitum}. (Dauben, p. 146)  \end{quote}

\section{Absolute truth}

     Cantor was sure that he had built his mathematical theory on a solid metaphysical foundation, theologically grounded. Consequently, he was convinced that his theory of infinite numbers
and his theory of real numbers were not only the only possible, but the only correct ones. He firmly rejected any other attempts. 

\begin{quote} My theory stands as firm as a rock, every arrow directed against it will return quickly to its archer. How do I know this? Because I have studied it from all sides for many years; because I have examined all objections which have been made against the infinite numbers; and above all, because I have followed all its roots, so to speak, to the first infallible cause of all created things. (Dauben 1990, p. 298). \end{quote}
     
His unshakable confidence in his theory was infectious, until today most mathematicians do not doubt it. This conviction had another consequence. As a responsible man who had uncovered the truth about the nature of infinity, he considered his duty to persuade the Church to accept it, and thus avoid the danger of a theological error that might result from accepting another theory of infinity.

\begin{quote} In any case, it is necessary to submit the question of the truth of the \emph{Transfinitum} to a serious examination, for were if the case that I am right in asserting the truth or possibility of the \emph{Transfinitum}, then (without doubt) there would be a sure danger of religious error in holding the opposite position, for: \emph{an error concerning the created is carried over into an erroneous doctrine of God} (Thomas Aquinas: Summa contra Gentiles 11,3). (Dauben 1990, p. 232). \end{quote}

Cantor was convinced that the knowledge of infinite numbers had been revealed to him by God, who guided his steps from pure mathematics to an interest in theology and philosophy so that he could improved a proper understanding of God and nature. 

\begin{quote}  But now I thank God, the all-wise and all-good, that He always denied me the fulfilment of this wish [for a position at university in Berlin or G\H ottingen], for He thereby has constrained me, through a deeper penetration into theology, to serve Him and His Holy Roman Catholic Church better than I have been able with my exclusive preoccupation with mathematics. (Dauben 1990, p. 147). \end{quote}

\section{Acceptance of set theory}

     Cantor continued to work hard at the systematic developing of his set theory. At the same time, he looked for a forum where mathematicians could freely present their new results and discuss them without fear of a prejudiced condemnation of a small elite of academics in Berlin. 

At that time, he devoted a considerable effort to reorganise the \emph{Section for Mathematics and Astronomy of the Society of German Scientists and Physicians}. The energy and enthusiasm with which Cantor set about this work bore fruit. A permanent professional \emph{Deutsche Mathematiker-Vereinung} (DMV) was established and Cantor was elected as a president. Its journal \emph{Mathematische Annalen} has been published until today. At the first meeting of DMV in 1891, Cantor read his paper that presented a new effective proof of the uncountability of real numbers, today called the \emph{Cantor diagonal method}.

This generalized theorem is more powerful than the previous one and has far-reaching consequences. It claims that a power-set, i.e. a set of all subsets of a given set, has a cardinality greater than that original set. If the original set had a cardinality $\kappa$ then its power-set, the set of all its subsets, has the cardinality $2^\kappa$. Consequently, there is an unlimited ascending sequence of infinite cardinalities. The Continuum Hypothesis was extended to the more general conjecture that the power-set of a set of a cardinality $\aleph_n$ would have the cardinality $\aleph_{n+1}$.
$$2^{\aleph_n} = \aleph_{n+1}.$$ 

     The following \emph{Contributions to the Foundation of the Theory of Transfinite Sets}, published successively in 1895 and 1897 is Cantor's last major mathematical publication. It is also the best known because it summarizes the purely mathematical core of the theory without philosophical considerations. Its spirit is foreshadowed  by the three initial quotations The first one is the famous Newton phrase: \emph{Hypotheses non fingo.} Like the other quotations, it testifies to Cantor's firm convictions that his theory of infinite numbers is a mathematical theory of permanent value, as natural and necessary as the theory of finite numbers. Its laws are fixed and are no random hypotheses.    

     After the publication of the \emph{Contributions} and their almost immediate translation into French and Italian, Cantor's theory became well known and spread among mathematicians all over the world. The theory of infinite numbers was appreciated and mathematicians were divided into two camps, those in favour and those against. Cantor no longer had to face opposition to his work alone but many energetic young mathematicians joined him. Although not all problems were solved, recognition by the mathematical community was achieved.

     At the first \emph{International Congress of Mathematicians} in Z\H urich in 1897, important applications of Cantor's set theory to mathematical analysis, topology and measure theory were pointed out. Just as enthusiastically as the theory had been rejected at first, it was now accepted and developed. Mathematicians appreciated the advantages it brought and gradually began to rebuild the foundations of mathematics on its basis. New disciplines emerged that could only be developed thanks to sets. At the second \emph{International Congress} in Paris in 1900, David Hilbert, the leading mathematician of the time, presented the twenty-three greatest problems in mathematics, the first one was to find a solution of the Continuum Hypothesis.

On Cantor's foundations, 20th century mathematicians gradually built a a sophisticated set theory based on infinite cardinal and ordinal numbers that is beautiful, challenging and complex. They have invented methods to create the whole hierarchy of much larger cardinal numbers and are exploring the relationships among them.  

\section{Infinitesimal quantities}

     It is somewhat ironic that infinitesimal quantities, whose uncertain existence was an important impetus for creation of a mathematical theory of actual infinity, cannot be directly defined in a Cantorian set theory.\footnote{The exception is non-standard analysis. Abraham Robinson around 1960 constructed a model in set theory using the ultrafilter that reformulates the calculus using a notion of infinitesimal numbers. However, it is not a direct description, an ultrafilter is a sophisticated mathematical structure.}

     Infinitesimal quantities were used by Newton and Leibniz from the end of 17th century to formulate the important differential and integral calculus, the infinitesimal calculus. In the mathematics built on the Cantor set theory, however, infinitesimals are not generally used. One actually works with potential infinity. For example, a continuity of a function f at a real number $c$ is defined as follows: For any $\epsilon > 0$ there exists $\delta > 0$ such that for every $x$ such that if $|x - c| < \delta$ then $|f(x) - f(c)| < \epsilon$.  For every arbitrarily small positive $\epsilon$, it is possible to find a positive $\delta$ that satisfies the given condition. If one had at disposal infinitesimal quantities traditionally denoted by $\Delta x$, one could just define that for each $\Delta x$ the difference $|f(c + \Delta x) - f(c)|$ is infinitesimal.

Infinitely large numbers look like a good basis for creating infinitely small numbers. Some mathematicians, e.g. Stolz, du Bois-Reymond, Vivanti, tried to introduce it during Cantor's lifetime, for example as inverse values of infinitely large numbers: $\frac{1}{\omega}, \frac{1}{\omega+1}, \dots \frac{1}{2\omega}, \dots $. Just this way would be very difficult. Both ordinal and cardinal have different arithmetical properties than finite numbers and cannot be simply divided or subtracted.  

However, Cantor firmly rejected any such ideas. He reserved the strongest words to infinitesimal quantities of du Bois-Reymond  and others. 
According to Cantor they have no meaning and should be thrown away as \enquote{nothing than paper numbers} because they do not satisfy the Archimedean axiom. Of course, it does not apply to transfinite numbers either. Somewhat paradoxically, many of Cantor's objections can be turned against his own infinitely large numbers as well. 
     
Cantor had yet another reason for opposing infinitesimal numbers. Accepting them would have even complicated his already rather complicated proof of the Continuum Hypothesis, which he had spent his life trying in vain to achieve.

\section{Consistent and inconsistent multiplicities}

     As is well known, some paradoxes appeared in set theory already during Cantor's lifetime. For example, if we consider the set of all sets, then its power-set has a greater cardinality, but at the same time it should be only its subset. The question also arose whether all ordinal numbers form a set. Since they are well-ordered they would have to be denoted by another ordinal number greater than all the others. The same is true for cardinal numbers.

     Cantor tried to solve this situation in a letter to Dedekind in 1899 by dividing sets into \emph{consistent} and \emph{inconsistent} sets. An inconsistent multiplicity [\emph{vielheit}] is such that \enquote{the assumption of the simultaneous existence of all elements leads to a contradiction} so that the quantity cannot be seen as a whole, as \enquote{one completed thing}. As an example, he gave the multiplicity of \enquote{totality everything thinkable}, but he could also have in mind the multiplicity of all sets. He also proved that the set of ordinal numbers and the set of cardinal numbers are inconsistent. 
 On the other hand, if the totality of the elements of a multiplicity can be considered without contradiction as \enquote{being together}, so that they can be gathered together into \enquote{one thing} then Cantor calls it a consistent multiplicity or a \enquote{set}. 

     The following letter, dated at Halle on 28 August 1899, is perhaps indicative of Cantor's uncertainty as to what is actually the case with the consistency of infinite sets.The question he asks here is not a marginal one.  The \enquote{simultaneous existence} of all the elements of an infinite set does not imply anything else than the existence of an actual infinite.
   
\begin{quote}  The question remains: how I know that the well-ordered multiplicities or sequences to which I ascribe the cardinal numbers 
$$\aleph_0, \aleph_1, \dots \aleph_{\omega_0}, \dots \aleph_{\omega_1}, \dots $$ are really \enquote{sets} in the sense of the word I have explained, i.e. \enquote{consistent multiplicities}. Is it not thinkable that these multiplicities are already \enquote{inconsistent}, and that the contradiction arising from the assumption of a \enquote{being together of all their elements} has simply not yet been made noticeable? My answer to this is that the same question can just as well be raised about finite multiplicities, and that a careful consideration will lead one to the conclusion that even for finite multiplicities no \enquote{proof} of their consistency is to be had. In other words, the fact of the \enquote{consistency} of finite multiplicities is a simple, unprovable truth; it is  (in the old sense of these words) \enquote{the axiom of finite arithmetic}. And in just the same way, the \enquote{consistency} of those multiplicities to which I attribute the alephs as cardinal numbers is \enquote{the axiom of extended transfinite arithmetic}. [Cantor, Letter to Dedekind, 28th August, 1899; GA 447-8] \end{quote}

 Cantor admits that what he has believed all his life, that is, consistency or \enquote{simultaneity} of all the elements of the transfinite numbers, he is unable to prove and has no choice but to establish as an axiom. He appeals to the fact that we tacitly assume the consistency of finite numbers in the same way.  At the same time it appears that Cantor viewed infinite sets as if they were finite. (Hallett ???). The distinction between finite and infinite is not admitted in the question of the guarantee of consistency.

\section{Axiomatic system} 

However, mathematicians\footnote{I write mathematicians, but it rather should be mathematicians of the mainstream. Some have not accepted Cantor's set theory, for example Henri Poincaré writes of it as an interesting pathological case and predicts that later generations will look upon it as a disease from which we have recovered. (Kline p. 1003).} who recognized great advantages that set theory offered were not about to give it up. Hilbert's famous quote is telling: \enquote{No one shall expel us from the paradise which Cantor has created for us.} (Kline, p. 1003). Bertrand Russell said Cantor’s work was \enquote{probably the greatest achievement of which our age can boast.} (Maor 1986, p. 23). Instead of rethinking the way in which actual infinity was introduced into set theory and the reasons that led to it, they restricted the notion of a set essentially on the basis of Cantor's definition of consistent sets.

The axiomatic systems on which of set theory was built were established. 
Two of them, Gödel-Bernays and Zermelo-Fraenkel, still the most widely used, are actually similar. They are intended to resolve old paradoxes and prevent new ones, and therefore they limit the meaning of the term set. Cantor's \enquote{inconsistent sets} are called classes, and it is not possible to perform certain operations on them (to determine their  ordinal and cardinal numbers, form their power-set, etc.). 

The axiomatic set theories no longer deal with the question of the existence of the actual infinite. Its positive solution is expressed axiomatically: \enquote{There is an infinite set}, or \enquote{The natural numbers form a set}. The Axiom of Infinity guarantees the character of the whole set theory established by Cantor. It expresses that there is an actual infinite set, a whole which we can treat as a single object.

Another fundamental axiom is the Power-Set Axiom. It guarantees that the set of all subsets of a given set is again a set. It expresses the belief that every set, even an infinite set, is transparent, and that all its parts, individually and all together, are sets. It implies that all real numbers form a set of greater cardinality than the natural numbers and that for any set we can form another set of a greater cardinality.

Axiom of Infinity guarantees the existence of an infinite increasing sequence of ordinal numbers, Axiom of Power-Set the existence of an infinite increasing sequence of cardinal numbers. The Continuum Hypothesis would connect these two sequences. However, this has not happened and cannot happen.

\section{Relative consistency}
 
Indeed, the Continuum Hypothesis, like many other important statements (e.g., the Axiom of Choice, the Axiom of Constructibility, the Diamond Principle), is an independent theorem. That is, it can neither be proved nor disproved from the axioms of set theory. One can consistently accept them or their negation.  

These independent theorems are not marginal; they are theorems important to properties of infinite sets and thus to the mathematics that is formulated in set theory. The situation is resolved by saying in advance which principles are accepted and then working with the theory in question. 

There is not a unique set theory, only its ground, which can be split into different branches, depending on one's choice. All that is certain is that one is working in a relatively consistent theory with respect to the the adopted axioms adopted. Mathematicians mostly work in Zermelo-Fraenkel set theory with the Axiom of Choice and they assume the Continuum Hypothesis. However, there is no clear and certain guideline to prioritise one branch over another.

Moreover, according to G\H odel's famous inconsistency theorem from 1931, if a mathematical theory that includes the axioms for natural numbers is consistent then it is not complete, that is, some of its valid statements are not provable, and vice versa, a complete theory is not consistent. 

Nobody can claim anymore, as Cantor did, absolute truth of the set theory, and no one does. This question is left aside. It is as if inconsistency were more important than truth, and as if truth lost its meaning, certainly had less meaning than consistent description of the world. This is, after all, one of the sources of postmodern relativism.

\section{Epilogue}

  I entirely agree with Cantor's words.    

\begin{quote} I believe that metaphysics and mathematics ought properly to be in relation to each other, and that at the time of their decisive achievements they were in mutual unity. Yet, unfortunately, as history shows, there is very soon a discord between them, which persists for generations, and which may grow to the point where the antagonistic brothers no longer know, or even want to know, that they are bound to each other.  (Cantor 1885) \end{quote}

      Philosophy and mathematics should be in harmony with each other, especially on such an important issue as infinity. Cantor based his philosophical belief on neo-Thomistic theology, which is not even a leading theological direction for the Catholic Church today. However, this had an indirect effect on the set theory he founded.
 
 Cantor, with his courageous work, gave mathematics the wings to take off in a new direction, but in his preoccupation with the battle for the actual infinite he gave it a greater range. Infinity somehow stiffened in Cantorian theory, multiplied and expanded in a strange way, and formed a marvellous but partly useless labyrinth of infinite cardinal and ordinal numbers.

The large cardinal numbers are a necessary consequence of set axioms. They are used only in purely theoretical mathematics and even rarely but they must be taken into account. However, they lost the correspondence with the physical world. They are not necessary for natural sciences, for which countability and continuum is sufficient. Both ordinal and cardinal numbers can only be increased: added, multiplied and powered, one can go ever higher. Hierarchy of large cardinals is thus rather \emph{l'art pour l'art}, beautiful but probably superfluous.

 Actual infinity belongs to mathematics, but God as its guarantor does not. If we reject it entirely, we considerably limit our research. But it is not necessary to accept them in Cantor's form. Set theory can also be built on another philosophical ground. 

A remarkable attempt to create an alternative set theory was made in Prague in the 1970s. (Vopěnka 1979). 
Its author, Petr Vopěnka, drew on phenomenological philosophy, various non-standard theories that appeared during the 20th century, and a critique of Cantor's set theory. He  interprets infinity on the basis of human evidence using the Husserlian concept of a horizon and  endeavours to preserve correspondence between the natural world and mathematical objects. 


\end{document}